\documentclass[11pt,reqno]{amsart}
\usepackage{amssymb,amscd,amsbsy}
\usepackage{amssymb,amscd,amsbsy,mathrsfs}
\newcommand{\lb}{\linebreak}

\renewcommand{\a}{\alpha}

\newcommand{\e}{\varepsilon}
\newcommand{\vk}{\varkappa}
\newcommand{\z}{\zeta}

\renewcommand{\l}{\lambda}

\newcommand{\s}{\sigma}
\renewcommand{\t}{\tau}

\newcommand{\f}{\varphi}
\renewcommand{\o}{\omega}

\newcommand{\D}{\Delta}

\renewcommand{\O}{\Omega}

\newcommand{\E}{{\mathscr E}}

\newcommand{\h}{{\mathscr H}}

\newcommand{\K}{{\mathscr K}}

\newcommand{\C}{{\Bbb C}}
\newcommand{\T}{{\Bbb T}}
\newcommand{\pp}{{\Bbb P}}
\newcommand{\dd}{{\Bbb D}}
\newcommand{\R}{{\Bbb R}}
\newcommand{\Z}{{\Bbb Z}}

\newcommand{\0}{{\boldsymbol{0}}}

\newcommand{\bs}{\boldsymbol}

\newcommand{\bS}{{\boldsymbol S}}

\newcommand{\rf}[1]{(\ref{#1})}

\newcommand{\df}{\stackrel{\mathrm{def}}{=}}

\newcommand{\Ker}{\operatorname{Ker}}
\newcommand{\re}{\operatorname{Re}}

\newcommand{\clos}{\operatorname{clos}}

\newcommand{\trace}{\operatorname{trace}}

\newcommand{\const}{\operatorname{const}}
\newcommand{\tr}{\operatorname{trace}}
\newcommand{\eeq}{\end{equation}}
\newcommand{\beq}{\begin{equation}}
\newcommand{\bay}{\begin{eqnarray}}
\newcommand{\ba}{\begin{align*}}
\newcommand{\ea}{\end{align*}}
\newcommand{\ey}{\end{eqnarray}}
\newcommand{\bey}{\begin{eqnarray*}}
\newcommand{\eey}{\end{eqnarray*}}

\newcommand{\be}{\infty}

\newcommand{\bl}{\blacksquare}

\newcommand{\Range}{\operatorname{Range}}

\newcommand{\Pf}{{\bf Proof. }}
\newcommand{\im}{\operatorname{Im}}
\renewcommand{\re}{\operatorname{Re}}
\newcommand{\ov}{\overline}

\newtheorem{thm}{\hspace{\parindent}Theorem}[section]

\newtheorem{cor}[thm]{\hspace{\parindent}Corollary}
\newtheorem{lem}[thm]{\hspace{\parindent}Lemma}

\theoremstyle{remark}

\newtheorem*{rem*}{Remark}

\newcommand\CA{{\rm C}_{\rm A}}
\newcommand{\OL}{{\rm OL}}
\newcommand{\OLA}{{\rm OL}_{\rm A}}

\newcommand{\ri}{{\rm i}}

\newcommand\mM{\mathcal{M}}

\newcommand\dg{\frak D}

\newcommand\cB{{\mathcal{B}}}

\newcommand\mB{\mathcal{B}}

\begin{document}

\newcommand{\vse}{\vspace{.2in}}
\numberwithin{equation}{section}

\title{Analytic operator Lipschitz functions in the disk and a trace formula for functions of contractions}

\author{M.M. Malamud, H. Neidhardt and V.V. Peller}

\maketitle

\footnotesize
\hfill To the memory of M.S. Agranovich, a remarkable mathematician and a remarkable personality
\normalsize

\begin{abstract}
In this paper we prove that for an arbitrary pair $\{T_1,T_0\}$ of contractions on Hilbert space with trace class difference, there exists a function $\bs\xi$ in $L^1(\T)$ (called a spectral shift function for the pair $\{T_1,T_0\}$ ) such that 
the trace formula $\trace(f(T_1)-f(T_0))=\int_\T f'(\z)\bs{\xi}(\z)\,d\z$
holds for an arbitrary operator Lipschitz function $f$ analytic in the unit disk.
\end{abstract}

\maketitle


\setcounter{section}{0}
\section{\bf Introduction}
\setcounter{equation}{0}
\label{In}

\medskip

The notion of spectral shift function was introduced by physicist I.M. Lifshits in \cite{L}. Later M.G. Krein elaborated the notion of spectral shift function in \cite{Kr} a most general situation; he showed that for self-adjoint operators $A_0$ and $A_1$ with trace class difference, there exists a unique real function 
$\bs{\xi}=\bs{\xi}_{A_1,A_0}$ in $L^1(\R)$ (it is called the {\it spectral shift function for the pair} $\{A_1,A_0\}$) such that the following trace formula holds:
\bay
\label{foslsso}
\trace\big(f(A_1)-f(A_0)\big)=\int_\R f'(t)\bs{\xi}_{A_1,A_0}(t)\,dt
\ey
for sufficiently nice functions $f$.

To prove the existence of the spectral shift function, he
introduced  the concept of {\it perturbation determinant} $\Delta_{A_1/A_0}$ and proved
the inversion formula
 $$
 \bs{\xi}_{A,A_0}(t) = \frac{1}{\pi}\lim_{y\downarrow 0}\im (\log(\Delta_{A_1/A_0}(t +\ri y)))
\quad \mbox{for a.e.} \quad  t \in \R,
$$
where  $\Delta_{A_1/A_0}(\z)\df\det(I + (A-A_0)(A_0-\z I)^{-1})$ (see \cite{Kr} and
\cite{Ya}).

Later Krein extended in \cite{K64} formula \eqref{foslsso} to the class $W_1(\R)$  of functions whose derivative is the Fourier transform of a complex Borel measure.

Krein also observed in \cite{Kr} that the right-hand side of \rf{foslsso} makes sense for arbitrary Lipschitz functions $f$ and posed the problem to describe the class of functions, for which formula \rf{foslsso} holds for all pairs of self-adjoint operators with difference of trace class $\bS_1$.

It turned out that trace formula \rf{foslsso} cannot be generalized to the class of all Lipschitz functions. Indeed, it was shown in \cite{F} that there exist a Lipschitz function $f$ on $\R$ and self-adjoint operators $A_1$ and $A_0$ such that $A_1-A_0\in\bS_1$, but $f(A_1)-f(A_0)\not\in\bS_1$. In \cite{Pe1} and \cite{Pe3} it was proved that \rf{foslsso} holds for functions $f$ in the (homogeneous) Besov space
$B_{\be,1}^1(\R)$ and does not hold unless $f$ locally belongs to the Besov space
$B_{1,1}^1$

Krein's problem was completely resolved recently in \cite{Pe6}. It was shown in \cite{Pe6} that trace formula \rf{foslsso} holds for arbitrary pairs $\{A_1,A_0\}$ of not necessarily bounded self-adjoint operators with trace class difference if and only if $f$ is an {\it operator Lipschitz function}, i.e., the inequality $$
\|f(A)-f(B)\|\le\const\|A-B\|
$$
holds for arbitrary self-adjoint operators $A$ and $B$.

In \cite{Kr2} Krein introduced the notion of spectral shift function for pairs of unitary operators with trace class difference. Namely, for a pair of unitary operators $\{U_1,U_0\}$ with trace class difference, he proved that there exists a real function $\bs{\xi}$ in $L^1(\T)$, unique modulo an additive constant, (called a {\it spectral shift function for} $\{U_1,U_0\}$ ) such that the trace formula 
\bay
\label{fosuo}
\trace\big(f(U_1)-f(U_0)\big)=\int_\T f'(\z)\bs{\xi}(\z)\,d\z
\ey
holds for functions for $f$ whose derivative $f'$ has absolutely convergent Fourier series.

In the recent paper \cite{AP+} an analog of the result of \cite{Pe6} was obtained: 
the class of functions $f$, for which formula \rf{fosuo} holds for arbitrary pairs $\{U_1,U_0\}$ of unitary operators with $U_1-U_0\in\bS_1$ coincides with the class $\OL_\T$ of operator Lipschitz functions on $\T$.

We are going to extend in this paper trace formulae \rf{foslsso} and \rf{fosuo}
to the case of pairs of maximal dissipative operators (m-dissipative) and pairs of contractions.


The Sz.-Nagy--Foia\c s functional calculus (see \cite{SNF})
associates with each function $f$ in the {\it disk-algebra} $\CA$ (i.e., the space of functions analytic in the disk $\dd$ and continuous in its closure) the operator $f(T)$. This functional calculus is linear and multiplicative and the von Neumann inequality 
$
\|f(T)\|\le\max\{|f(\z)|:~\z\in\C,~|\z|\le1\}
$
holds for $f\in\CA$.

When we proceed to the case of contractions, it is natural to divide the problem in  two
problems: (i) the existence of a spectral shift function and a trace formula for resolvents; (ii) a
description of the largest class  of functions for which the trace formula 
holds.

First generalizations  of formula \eqref{foslsso} to the case of pairs $\{A_1, A_0\}$ with
an m-accumul\-ative (dissipative) operator $A_1$ were obtained by Rybkin \cite{Ryb84} and
Krein \cite{Kr87}.
For
instance, Krein treated the case when $A_0 = A^*_0$,  $A_1 = A_0 - {\rm i}V$,
$V\ge 0$ and $V\in \bS_1$ and proved in \cite{Kr87} an analog of formula \eqref{foslsso} with right-hand side $\int_\R f'(t)\,d\nu(t)$ for
 a complex Borel measure $\nu$ and for functions $f$ of class $W_1^+$, i.e., functions $f$ whose derivative is the Fourier transform of a complex measure supported in $[0,\be)$.
The most complete result on the existence of a spectral shift function for  pairs $\{L_1,L_0\}$ of m-accumulative resolvent
comparable operators satisfying $\rho(L_0)\cap \C_+\not = \varnothing$
was obtained in \cite{MalNei2015}. Under these
assumptions, formula \eqref{foslsso} was proved  for resolvents.  In the case of
additive $\bS_1$ perturbations  trace formula \eqref{foslsso} was proved for the class
$W_1^+$ (see \cite{MalNei2015}).

The purpose of this paper is to prove that for a pair $\{T_1,T_0\}$  of contractions satisfying $T_0-T_1\in\bS_1$, there exists a function $\bs\xi$ in $L^1(\T)$ (called a {\it spectral shift function for the pair} $\{T_1,T_0\}$ ) such that the following trace formula
\bay
\label{tfszha}
\trace\big(f(T_1)-f(T_0)\big)=\int_\T f'(\z)\bs{\xi}(\z)\,d\z
\ey
holds for an arbitrary operator Lipschitz function $f$ analytic in $\dd$.

To obtain the main result, we combine two approaches. The first approach is based on double operator integrals with respect to semi-spectral measures and uses an idea of \cite{BS2}. It leads to a trace formula 
$\trace\big(f(T)-f(R)\big)=\int_\T f'(\z)\,d\nu(\z)$ for an arbitrary operator Lipschitz function $f$ analytic in $\dd$, where $\nu$ is a Borel measure on $\T$.

The second approach develops ideas of Krein \cite{Kr} - \cite{Kr87} and relies on the study of the perturbation determinant for pairs
$\{L_1,L_0\}$ of m-accumulative operators and allows us to prove the trace formula for
resolvents without  additional restrictions on $\{L_1,L_0\}$ imposed in \cite{MalNei2015}.

We denote by $\mB(\h)$ the set of bounded linear operators on a Hilbert space $\h$ and by $\bS_p$ the
Schatten-von Neumann ideal. In particular,  $\bS_1$ is the trace
class.  For a closed densely
defined operator $Q$, we use the notation $\rho(Q)$, $\s_{\rm p}(Q)$,  and 
$\s_{\rm c}(Q)$
denote the resolvent set and the point and the continuous spectra  of $Q$. Recall
that 
$\s_{\rm c}(Q) =\{\l\notin\s_{\rm p}(Q):~\Range(Q-\l I)\ne\clos\Range(Q-\l I) =
\h$\}.

\section{\bf Double operator integrals and Schur multipliers}
\setcounter{equation}{0}
\label{dois}

\medskip

Double operator integrals
$$
\iint\Phi(x,y)\,dE_1(x)Q\,dE_2(y)
$$
were introduced by Yu.L. Daletskii and S.G. Krein in \cite{DK}. Later
Birman and Solomyak elaborated their beautiful theory of double oparator integrals \cite{BS} (see also \cite{AP} and references there). Here $\Phi$ is a bounded measurable function, $E_1$ and $E_2$ are spectral measures on Hilbert space, and
$Q$ is a bounded linear operator. Such double operator integrals are defined for arbitrary bounded measurable functions $\Phi$ if $Q$ is a Hilbert--Schmidt operator. If $Q$ is an arbitrary bounded operator, then for the double operator integral to make sense, $\Phi$ has to be a Schur multiplier with respect to $E_1$ and $E_2$, (see \cite{Pe1} and \cite{AP}). It is well known (see \cite{Pe1} and \cite{AP}) that $\Phi$ is a Schur multiplier if and only if it admits a representation
$$
\Phi(x,y)=\sum_{n\ge0}\f_n(x)\psi_n(y),
$$
where the $\f_n$ and $\psi_n$ are measurable functions such that
$$
\sum_{n\ge0}|\f_n|^2\in L^\be_{E_1}\quad\mbox{and}\quad
\sum_{n\ge0}|\psi_n|^2\in L^\be_{E_2}.
$$

In this paper we deal with double operator integrals with respect to {\it semi-spectral measures}
$$
\iint\Phi(x,y)\,d\E_1(x)Q\,d\E_2(y).
$$
Such double operator integrals were introduced in \cite{Pe3} (see also \cite{Pe6}).
We refer the reader to the recent survey \cite{AP} for detailed information.

If $T$ is a contraction on a Hilbert space $\h$, it has a {\it minimal unitary dilation} $U$, i.e., $U$ is a unitary operator on a Hilbert space $\K$, $\K\supset\h$, $T^n=P_\h U^n\big|\h$ for $n\ge0$ and $\K$ is the closed linear span of
$U^n\h$, $n\in\Z$ (see \cite{SNF}). Here $P_\h$ is the orthogonal projection onto $\h$.
The {\it semi-spectral measure $\E_T$ of} $T$ is defined by
$$
\E_T(\D)\df P_\h E_U(\D)\big|\h,
$$
where $E_U$ is the spectral measure of $U$ and $\D$ is a Borel subset of $\T$. It is easy to see that 
$
T^n=\int_\T\z^n\,d\E_T(\z),~ n\ge0.
$

\section{\bf The role of divided differences}
\setcounter{equation}{0}
\label{Haage}

\medskip

A function $f$ in the disk-algebra $\CA$ is called {\it operator Lipschitz} if
$$
\|f(T)-f(R)\|\le\const\|T-R\|
$$
for all contractions $T$ and $R$. We use the notation 
$\OLA$ for the class of operator Lipschitz functions in $\CA$. 
For $f\in\OLA$, we denote by $\dg f$ the divided difference defined by
$$
(\dg f)(\z,\t)=\dfrac{f(\z)-f(\t)}{\z-\t}\quad\mbox{for}\quad\z\ne\t\quad\mbox{and}\quad
(\dg f)(\z,\t)=f'(\z)\quad\mbox{for}\quad\z=\t.
$$
Recall that by a theorem of Johnson and Williams \cite{JW}, operator Lipschitz functions must be differentiable everywhere on $\T$.

We are going to use the following representation of the divided difference 
$\dg f$ for functions in $\OLA$ (see \cite{AP}, Theorems 3.9.1 and 3.9.2):

\medskip

{\it Let $f$ be a function analytic in $\dd$. Then $f\in\OLA$ if and only if
$f$ admits a representation 
\bay
\label{Haagerpre}
(\dg f)(z,w)=\sum_{n\ge1}\f_n(z)\psi_n(w),\quad z,\;w\in\dd,
\ey
where $\f_n$ and $\psi_n$ are functions in $\CA$  such that
\bay
\label{proisumm}
\Big(\sup_{z\in\dd}\sum_{n\ge1}|\f_n(z)|^2\Big)\Big(\sup_{w\in\dd}\sum_{n\ge1}|\psi_n(w)|^2\Big)<\be
\ey
If $f\in\OLA$, then the functions $\f_n$ and $\psi_n$ can be chosen so that
the left-hand side of {\em\rf{proisumm}} is equal to $\|f\|_{\OLA}$.}

Let $T_0$ and $T_1$ be contractions on Hilbert space, and let $\E_0$ and $\E_1$ be their semi-spectral measures. Suppose now that $f\in\OLA$. 
Consider a representation of $\dg f$ in the form \rf{Haagerpre}, where $\f_n$ and $\psi_n$ are functions in $\CA$ satisfying \rf{proisumm}. Then for a bounded linear operator $K$, we have
\bay
\label{dvoopiryad}
\iint_{\T\times\T}\big(\dg f)(\z,\t)\,d\E_1(\z)K\,d\E_0(\t)
=\sum_{n=1}^\be\f_n(T_1)K\psi_n(T_0)
\ey
(see Section 3.9 of \cite{AP}).
This implies (see Theorem 3.9.9 of \cite{AP}) that
\bay
\label{razdvoopi}
f(T_1)-f(T_0)=
\iint_{\T\times\T}\big(\dg f)(\z,\t)\,d\E_1(\z)(T_1-T_0)\,d\E_0(\t).
\ey

\section{\bf Differentiation in the strong operator topology}
\setcounter{equation}{0}
\label{difvsil}

\medskip

The following theorem is a generalization of Theorem 3.5.6 of \cite{AP} and
Theorem 3.1 of \cite{AP+}.

\begin{thm}
\label{operdifszha}
Let $f\in\OLA$ and let $T_0$ and $T_1$ be contractions on Hilbert space and
$T_t\df T_0+t(T_1-T_0)$, $0\le t\le1$. Then
$$
\lim_{s\to0}\frac1s\big(f(T_{t+s})-f(T_t)\big)
=\iint_{\T\times\T}\big(\dg f)(\z,\t)\,d\E_t(\z)(T_1-T_0)\,d\E_t(\t)
$$
in the strong operator topology,
where $\E_t$ is the semi-spectral measure of $T_t$.
\end{thm}

We need the following elementary fact:

\medskip

{\bf Lemma 3.5.9 of \cite{AP}.} {\it Let $\{X_n\}_{n\ge1}$ be a sequence of bounded linear operators on a Hilbert space $\h$ and let $\{u_n\}_{n\ge1}$ be a sequence of vectors in $\h$. Suppose that
$$
\sum_{n\ge1}X_nX_n^*\le a^2I\quad\mbox{and}\quad\sum_{n\ge1}\|u_n\|^2\le b^2
$$
for  $a,\,b\ge0$. Then the series $\sum_{n\ge1}X_nu_n$ converges weakly and}
$
\Big\|\sum_{n\ge0}X_nu_n\Big\|\le ab.
$

\medskip

{\bf Proof of Theorem \ref{operdifszha}.} 
Let $\f_n$ and $\psi_n$ are functions in $\CA$ that satisfy \rf{Haagerpre} and
\rf{proisumm}. It follows from \rf{razdvoopi} and \rf{dvoopiryad} that it suffices to prove that
$$
\lim_{s\to0}\sum_{n\ge1}\f_n(T_{t+s})(T_1-T_0)\psi_n(T_t)=
\sum_{n\ge1}\f_n(T_{t})(T_1-T_0)\psi_n(T_t).
$$
Note that the limit should be taken in the strong operator topology while the series converge in the weak operator topology. 

Let $u$ be a vector of the Hilbert space. We have to prove that
$$
\lim_{s\to0}\big(\f_n(T_{t+s})-\f_n(T_t)\big)(T_1-T_0)\psi_n(T_t)u=\0
$$
in the Hilbert space norm.

We may assume that $\|u\|=1$, $\|f\|_{\OLA}\le1$, $\sum_{n\ge1}|\f_n|^2\le1$ and
$\sum_{n\ge1}|\psi_n|^2\le1$.

Put $u_n\df(T_1-T_0)\psi_n(T_t)u$. We have 
$$
\sum_{n\ge1}\|u_n\|^2\le\|T_1-T_0\|^2\sum_{n\ge1}\|\psi_n(T_t)u\|^2.
$$
Suppose now that $U$ is a unitary dilation of $T_t$. We have
$$
\sum_{n\ge1}\|\psi_n(T_t)u\|^2\le\sum_{n\ge1}\|\psi_n(U)u\|^2
=\sum_{n\ge1}(|\psi_n|^2(U)u,u)^2\le1,
$$
and so
$
\sum_{n\ge1}\|u_n\|^2\le\|T_1-T_0\|^2.
$
Let $\e>0$. We can select a positive integer $N$ such that 
$
\sum_{n>N}\|u_n\|^2<\e^2.
$
Let us show that
\bay
\label{dvanera}
\sum_{n\ge1}\f_n(T_t)\big(\f_n(T_t)\big)^*\le I\quad\mbox{and}\quad
\sum_{n\ge1}\f_n(T_{t+s})\big(\f_n(T_{t+s})\big)^*\le I.
\ey
Indeed, let $U$ be a unitary dilation of $T_t$ on a Hilbert space $\K$ that contains $\h$. It is easy to verify that 
$\big\|\big(\f_n(T_t)\big)^*v\big\|_\h\le\big\|\big(\f_n(U)\big)^*v\big\|_\K$.
We have
$$
\sum_{n\ge1}\big(\f_n(T_t)\big(\f_n(T_t)\big)^*v,v\big)
=\sum_{n\ge1}\big\|\big(\f_n(T_t)\big)^*v\big\|^2_\h
\le\sum_{n\ge1}\big\|\big(\f_n(U)\big)^*v\big\|^2_\K
\le\|v\|^2_\h
$$
because $\sum_{n\ge1}\|\f_n\|^2\le1$. This means that
$
\sum_{n\ge1}\f_n(T_t)\big(\f_n(T_t)\big)^*\le I.
$
In the same way we can prove that
$
\sum_{n\ge1}\f_n(T_{t+s})\big(\f_n(T_{t+s})\big)^*\le I.
$

It follows from the above lemma that for all $t$,
$$
\Big\|\sum_{n\ge1}\big(\f_n(T_{t+s})-\f_n(T_t)\big)u_n
\Big\|\le2\e.
$$

It is easy to see that
$
\|\f_n(T_{t+s})-\f_n(T_t)\|\to0\quad\mbox{as}\quad s\to0
$
(it suffices to approximate $\f_n$ be analytic polynomials). Thus
$$
\Big\|\sum_{n=1}^N\big(\f_n(T_{t+s})-\f_n(T_t)\big)u_n\Big\|
\le\|T_1-T_0\|\sum_{n=1}^N\|\f_n(T_{t+s})-\f_n(T_t)\|<\e
$$
if $s$ is sufficiently close to $0$. Thus
$
\Big\|\sum_{n\ge1}\big(\f_n(T_{t+s})-\f_n(T_t)\big)u_n
\Big\|<3\e
$
for $s$ sufficiently close to $0$ which completes the proof. $\bl$

\section{\bf A trace formula for double operator integrals}
\setcounter{equation}{0}
\label{forsledvo}

\medskip

\begin{thm}
\label{sleddvoopi}
Let $f\in\OLA$. Suppose that $T$ is a contraction with semi-spectral measure $\E$ and $K$ is a trace class operator  on Hilbert space. Then
$$
\trace\Big(\,\,\iint_{\T\times\T}\big(\dg f\big)(\z,\t)\,d\E(\z)K\,d\E(\t)\Big)=
\int_\T f'(\z)\,d\mu(\z),
$$
where $\mu$ is a complex Borel measure on $\T$ defined by
$
\mu(\D)=\trace(K\E(\D)).
$
\end{thm}

\Pf Let $\f_n$ and $\psi_n$ be functions in $\CA$ satisfying \rf{Haagerpre} and \rf{proisumm}. By \rf{dvoopiryad}, 
$$
\iint_{\T\times\T}\big(\dg f\big)(\z,\t)\,d\E(\z)K\,d\E(\t)=
\sum_{n\ge1}\f_n(T)K\psi_n(T),
$$
and so
\begin{align*}
\trace&\left(\,\,\iint_{\T\times\T}\big(\dg f\big)(\z,\t)\,d\E(\z)K\,d\E(\t)\right)=
\sum_{n\ge1}\trace\big(\f_n(T)K\psi_n(T)\big)\\[.2cm]
&=\sum_{n\ge1}\trace\big(K(\f_n\psi_n)(T)\big)
=\trace\Big(K\int_\T\Big(\sum_{n\ge1}\f_n\psi_n\Big)(\z)\,d\E(\z)
\Big)\\[.2cm]
&=\trace\left(K\int_\T\big(\dg f\big)(\z,\z)\,d\E(\z)\right)
=\trace\left(K\int_\T f'(\z)\,d\E(\z)\right)
=\int_\T f'(\z)\,d\mu(\z).
\end{align*}
This completes the proof. $\bl$

\section{\bf A trace formula for functions of contractions}
\setcounter{equation}{0}
\label{formszhat}

\medskip

\begin{thm}
\label{fsdm}
Let $T_0$ and $T_1$ be contractions on Hilbert space such that $T_1-T_0\in\bS_1$.
Then there exists a complex Borel measure $\nu$ on $\T$ such that the trace formula
\bay
\label{formsledlyaszha}
\trace\big(f(T_1)-f(T_0)\big)=\int_\T f'(\z)\,d\nu(\z)
\ey
holds for every $f$ in $\OLA$.
\end{thm}

\Pf Let $K\df T_1-T_0$. Consider the family of contractions $T_t\df T+tK$, $0\le t\le1$.
By Theorem \ref{operdifszha}, the function $t\mapsto f(T_t)$ is differentiable in the strong operator topology. Put
$$
Q_t\df\lim_{s\to0}\frac1s\iint_{\T\times\T}
\big(\dg f)(\z,\t)\,d\E_t(\z)K\,d\E_t(\t).
$$
Since $K\in\bS_1$ and $f\in\OLA$, we have
$$
Q_t\in\bS_1,\;0\le t\le1,\quad\mbox{and}\quad\sup_{t\in[0,1]}\|Q_t\|_{\bS_1}<\be.
$$

It follows from the definition of $Q_t$ that the function
$t\mapsto Q_tu$ is measurable for every vector $u$. Then 
the function $t\mapsto\trace (Q_tW)$ is measurable for an arbitrary bounded operator $W$. In other words the $\bS_1$-valued function $t\mapsto Q_t$ is
weakly measurable. Then it must be strongly measurable because $\bS_1$ is separable (see \cite{Y}, Ch. V, \S\:4) and
$$
f(T_1)-f(T_0)=\int_0^1Q_t\,dt,
$$
where the integral is understood in the sense of Bochner.

By Theorem \ref{sleddvoopi}, 
$$
\trace Q_t=\int_\T f'(\z)\,d\nu_t(\z),
$$
where $\nu_t$ is the complex Borel measure on $\T$ defined by
$
\nu_t(\D)\df\trace\big(K\E_t(\D)\big)
$
for a Borel subset $\D$ of $\T$.
It would be tempting to define the measure $\nu$ by
$$
\nu(\D)=\int_0^1\nu_t(\D)\,dt
$$
and conclude that formula \rf{formsledlyaszha} holds for this measure $\nu$.
However, it is not quite clear why the function $t\mapsto\nu_t(\D)$ is measurable.

To overcome this problem, we consider the dual space to the disk-algebra $\CA$. It can naturally be identified with the quotient space $\mM/H^1$,
where $\mM$ is the space of complex Borel measures on $\T$ and $H^1$ is the Hardy class.
Consider the function $t\mapsto \dot\nu_t$, where $\dot\nu_t$ is the coset in 
$\mM/H^1$ that corresponds to $\nu_t$.
Let us show that the map $t\mapsto\dot\nu_t$ is continuous in the weak-* topology on $\mM/H^1$. Indeed, let $h\in\CA$. We have
$
\langle h,\dot\nu_t\rangle=\trace\big(Kh(T_t)\big).
$
The result follows from the fact that the map $t\mapsto h(T_t)$ is continuous in the operator norm. This can be proved easily by approximating $h$ by analytic polynomials. 

Consider now the integral
$
\int_0^1\dot\nu_t\,dt.
$
It is an element of $\mM/H^1$, and so it is equal to $\dot\nu$ for some $\nu\in\mM$. It remains to observe that $\nu$ satisfies equality \rf{formsledlyaszha}. $\bl$

\newcommand{\red}{\color{red}}

\newcommand\eps{\varepsilon}
\newcommand\ff{\varphi}
\newcommand\kk{\varkappa}
\newcommand\tg{\tan}
\newcommand\ctg{\cot}
\newcommand\equ{{\Longleftrightarrow}}
\newcommand\then{\Longrightarrow}
\newcommand\eset{\varnothing}

\def\RE{{\rm Re\,}}
\def\IM{{\rm Im\,}}
\def\ran{{\rm ran\,}}
\def\tr{{\rm tr\,}}
\def\oD{{\bC\setminus\overline{\bD}}}
\def\dom{{\rm dom}}


\newcommand{\bR}{{\mathbb{R}}}
\newcommand{\bC}{{\mathbb{C}}}
\newcommand{\bD}{{\mathbb{D}}}
\newcommand{\bP}{{\mathbb{P}}}
\newcommand{\bT}{{\mathbb{T}}}
\newcommand{\bN}{{\mathbb{N}}}
\newcommand{\bZ}{{\mathbb{Z}}}

\newcommand\gotD{\mathfrak{D}}
\newcommand\gotH{\mathfrak{H}}
\newcommand\gotK{\mathfrak{K}}
\newcommand\goth{\mathfrak{h}}
\newcommand\gotL{\mathfrak{L}}
\newcommand\gotN{\gN}
\newcommand\gotR{\gR}
\newcommand\gotS{\mathfrak{S}}
\newcommand\gotX{\mathfrak{X}}
\newcommand\gotY{\mathfrak{Y}}
\newcommand\gotZ{\mathfrak{Z}}
\newcommand\gotz{\mathfrak{z}}
\newcommand\gotm{\mathfrak{m}}

\newcommand{\ga}{{\alpha}}
\newcommand{\gd}{{\delta}}
\newcommand{\gD}{{\Delta}}
\newcommand{\gga}{{\gamma}}
\newcommand{\gG}{{\Gamma}}
\newcommand{\gF}{{\Phi}}
\newcommand{\gf}{{\phi}}
\newcommand{\gk}{{\kappa}}
\newcommand{\gK}{{\Kappa}}
\newcommand{\gl}{{\lambda}}
\newcommand{\gL}{{\Lambda}}
\newcommand{\gO}{{\Omega}}
\newcommand{\go}{{\omega}}
\newcommand{\gs}{{\sigma}}
\newcommand\gS{{\Sigma}}
\newcommand{\gth}{{\theta}}
\newcommand{\gT}{{\Theta}}

\newcommand\cA{{\mathcal{A}}}

\newcommand\cC{{\mathcal{C}}}
\newcommand\cD{{\mathcal{D}}}
\newcommand\cF{{\mathcal{F}}}
\newcommand\cH{{\mathcal{H}}}
\newcommand\cK{{\mathcal{K}}}
\newcommand\cN{{\mathcal{N}}}
\newcommand\cO{{\mathcal{O}}}
\newcommand\cP{{\mathcal{P}}}
\newcommand\cT{{\mathcal{T}}}
\newcommand\cV{{\mathcal{V}}}

\newcommand\tA{\widetilde{A}}
\newcommand\tS{\widetilde{S}}
\newcommand\tT{\widetilde{T}}

\newtheorem{proposition}[thm]{Proposition}
\newtheorem{corollary}[thm]{Corollary}
\newtheorem{lemma}[thm]{Lemma}
\newtheorem{definition}[thm]{Definition}
\newtheorem{example}[thm]{Example}
\newtheorem{remark}[thm]{Remark}

\newcommand{\slim}{\,\mbox{\rm s-}\hspace{-2pt} \lim}
\newcommand{\wlim}{\,\mbox{\rm w-}\hspace{-2pt} \lim}
\newcommand{\olim}{\,\mbox{\rm o-}\hspace{-2pt} \lim}

\newcommand\gH{{\mathfrak{H}}}
\newcommand\fS{{\mathfrak{S}}}

\def\wt#1{{{\widetilde #1} }}
\def\wh#1{{{\,\widehat #1\,} }}
\def\Ext{{\rm Ext\,}}





\section{\bf The existence of spectral shift functions  for a pair of contractions}

\medskip

{\bf Definition.}
A densely defined  operator $L$ in $\h$ is called {\it dissipative} if $\IM(Lx,x) \ge 0$ for $x$ in the domain
$\dom(L)$ of $L$. It is called {\it maximal dissipative} (or {\it m-dissipative)} if $L$ has no proper dissipative extension. An operator $L$ is called accumulative
(m-accumulative) if $-L$ is dissipative (m-dissipative).

\medskip

First, we present an improvement of Theorem 3.14 of \cite{MalNei2015} for pairs
$\{L_1,L_0\}$ of m-dissipative operators.

\begin{thm}
\label{thm:2.1}
Let $L_1$ and $L_0$ be m-dissipative operators  such that  $(\rho(L_0) \cup
\sigma_{\rm c}(L_0))\cap \bC_+ \not= \varnothing$. Suppose that $(L_1 +\ri I)^{-1} - (L_0 + \ri I)^{-1}\in
\bS_1$, i.e., $L_1$ and $L_0$ are resolvent comparable.   Then
there exists a complex-valued function $\bs\o \in L^1(\bR,(1+t^2)^{-1})$ (a spectral shift function for the pair $\{L_1,L_0\}$ ) such
that the following  trace formula holds 
\begin{equation}
\label{eq:2.1}
\trace\big((L_1-\t I)^{-1} - (L_0-\t I)^{-1}\big) = -
\int_\bR\frac{\bs\o(t)}{(t-\t)^2}dt,
\quad\ \t\in \bC_-.
\end{equation}
If  $L_0 = L_0^*$  (respectively, $L_1 = L_1^*$), then one can select $\bs\o$ such that  $\IM(\bs\o(t)) \ge 0$ $($respectively, $\IM(\bs\o(t)) \le 0)$  a.e. on  $\bR$.
  \end{thm}

The proof of Theorem \ref{thm:2.1} is given in the Appendix. In what follows we use the notation $\mathscr S\{L_1,L_0\}$ for the set of spectral shift functions
for $\{L_1,L_0\}$.

\begin{thm}
\label{main_theor}
Let $T_1$ and $T_0$ be contractions with $T_1- T_0\in \bS_1$. Then
there is a complex-valued function $\bs\xi$ in $L^1(\T)$ 
(a spectral shift function for  $\{T_1,T_0\}$ )
such that 
   \begin{equation}\label{eq:2.2}
\trace\big((T_1 - \l I)^{-1} - (T_0 - \l I)^{-1}\big) = - 
\int_\bT
\frac{\bs\xi(\z)}{(\z - \l)^2}d\z, \qquad  |\l|>1.
  \end{equation}
If  $T_0$  (respectively, $T_1$) is unitary, we can select $\bs\xi$ such that  $\IM(\bs\xi(\z)) \le 0$  $($respectively, $\IM(\bs\xi(\z)) \ge 0)$  for
almost all $\z \in \bT$.
  \end{thm}

\Pf (i). Let us first assume that  $\Ker(T_1+ I)  = \Ker(T_0 + I) = \{0\}$ and  $0 \in \rho(T_0) \cup \gs_c(T_0)$. Note that for a
contraction $T$, the condition $\Ker(T + I) = \{0\}$ is equivalent to $\Ker(T^* + I) =
\{0\}$.   Therefore the operators
   \begin{equation}
   \label{2.4_oper_H_and_H_0}
L_1= -\ri I + 2\ri(I+T_1)^{-1} \quad \mbox{and} \quad L_0= -\ri I + 2\ri(I+T_0)^{-1}
\end{equation}
exist and are densely defined.   One can easily verify
that  $L_1$ and $L_0$ are dissipative. Since the operators $L_1+\ri I$ and 
$L_0+\ri I$ are onto, they  are m-dissipative. Let us show that $L_1$ and $L_0$ satisfy the
hypotheses of Theorem   \ref{thm:2.1}.

Clearly, the map $\l \mapsto\t \df \ri(1-\l)(1+\l)^{-1}$  maps $\bD$ onto  $\bC_+$ and maps $\bC\setminus\clos\dd$  onto 
$\bC_-$.
In accordance with  \eqref{2.4_oper_H_and_H_0},
  \begin{equation}\label{eq:2.5}
    \begin{split}
(L_1-\t I)^{-1} &= -\frac{1+\l}{2\ri}I - \frac{(1+\l)^2}{2\ri}(T_1-\l I)^{-1},\\
(L_0 -\t I)^{-1} &= -\frac{1+\l}{2\ri}I - \frac{(1+\l)^2}{2\ri}(T_0 -\l I)^{-1},
\end{split}
\quad |\l|>1 
   \end{equation}
where $\t\in \bC_-$.
It follows that $L_1$ and $L_0$ are resolvent comparable.

Let us show that $\ri \in \rho(L_0) \cup \s_{\rm c}(L_0)$. By \eqref{2.4_oper_H_and_H_0},
$T_0 = (\ri I-L_0)(\ri I+L_0)^{-1}$, and so the condition $\Ker T_0 = \{\0\}$  is
equivalent to $\Ker(L_0-\ri I) = \{\0\}$.  Similarly, $T_0$ has dense range
if and only if $L_0-\ri I$ has dense range.
This proves that
$\ri\in \rho(L_0)\cup\s_{\rm c}(L_0)$. By Theorem \ref{thm:2.1}, there exists
$\bs\o \in L^1(\bR,(1+t^2)^{-1})$ such that trace formula
\eqref{eq:2.1} holds.  Combining this formula with \eqref{eq:2.5} and setting
\begin{equation*}
t\df\ri(1-\z)(1+\z)^{-1} =\ri(1-e^{i\theta})(1+e^{i\theta})^{-1}\in\R,  \qquad \z = e^{i\theta} \in \bT, \quad\theta \in (-\pi,\pi],
\end{equation*}
and $\bs\xi(\z)\df\bs\o\left(\ri(1-\z)(1+\z)^{-1}\right)$,  i.e.,
\begin{equation}
\label{2.6_connection_of_two_SSF}
\bs\xi(e^{i\theta})=\bs\o\left(\ri(1- e^{i\theta})(1+e^{\ri\theta})^{-1}\right) =
\bs\o(\tg(\theta/2)), \qquad \theta\in  (-\pi,\pi],
\end{equation}
we obtain
\begin{align}
-\frac{(1+\l)^2}{2i}&\trace\left((T_1-\l I)^{-1}-(T_0-\l I)^{-1}\right) =
\trace((L_1-\t I)^{-1} - (L_0-\t I)^{-1})  \nonumber  \\
 &= -\int_\bR\frac{\bs\o(t)}{(t-\t)^2}dt = \frac{(1+\l)^2}{2\ri}
\int_\bT \frac{\bs\xi(\z)}{(\z-\l)^2}d\z, \qquad  |\l|>1.
\end{align}
This proves \eqref{eq:2.2}. Moreover,  
$\int_{-\pi}^{\pi}|\bs\xi(e^{\ri\theta})|\,d\theta = 2\int_{\R}
|\bs\o(t)|(1+t^2)^{-1}dt < \infty$, i.e., $\bs\xi\in L^1(\T)$.

(ii). Assume now that at least one of the conditions  $\Ker(T_1-I) =
 \{\0\}$ and \lb$\Ker(T_0-I)=\{\0\}$ is violated.
We construct an auxiliary pair  $\{\widehat T_1, \widehat T_0\}$
of contractions  with $\Ker(\widehat T_1-I)=\Ker(\widehat T_0-I)=\{\0\}$.

First, we observe that for each eigenvalue $\z\in\sigma_{\rm p}(T)$, 
$|\z| =1$,
the eigensubspace $\Ker(T-\z I)$ is also the eigensubspace
of the minimal unitary dilation $U$, i.e., $\Ker(U-\z I) = \Ker(T-\z I)$ (see
\cite{SNF}). Therefore the set $\sigma_{\rm p}(T)\cap\bT$ is at most countable. This
implies the existence of a number $\varkappa \in \bT$ such that the contractions
$\widehat T \df \varkappa T$ and $\widehat T_0 \df \varkappa T_0$ have the
required properties
 $$
\Ker(\widehat T+ I) = \Ker(\widehat T^*+I) = \{0\} \qquad \text{and} \qquad
\Ker(\widehat T_0 + I) = \Ker(\widehat T^*_0 + I) = \{0\}.
$$
In accordance with (i),  there exists a function $\widehat{\bs\xi}\in
L^1(\bT)$ such that
\begin{equation}
\label{2.8_trace}
\trace\left((\widehat T -\hat\l I)^{-1} - (\widehat T_0-\hat\l I)^{-1}\right) = -
\int_\T \widehat{\bs\xi}(\hat\z)(\hat\z -\hat\l)^{-2}d\hat\z, \qquad |\hat\l|>1.
\end{equation}
Setting $\l= \overline{\vk}\hat\l$, $\z=\overline{\vk}\hat\z$ and $\bs\xi(\z)= \widehat{\bs\xi}(\vk\z)$,
we obtain \eqref{eq:2.2}.

(iii). Suppose now that the condition  $0\in\rho(T_0) \cup\s_{\rm c}(T_0)$ is violated,
i.e., $0\in\s_{\rm p}(T_0)\cup \s_{\rm p}(T_0^*)$. Let us first assume that
$0\in\s_{\rm p}(T_0)$. Consider the Hilbert space $\h\oplus\K$, where $\h$ is the
space, on which $T_0$ and $T_1$ are defined and $\K$ is an infinite-dimensional Hilbert
space. Consider the operators
\begin{equation}
\label{2.10}
\check T_1= T_1\oplus\0 \quad \mbox{and} \quad \check T_0 = T_0 \oplus \0.
\end{equation}
on this space.
Clearly, $\check T_1-\check T_0\in\bS_1$. Let
 $T_0 = V_0|T_0|$ be the polar decompostion of $T_0$, where $V_0$ is a partial isometry that maps
$\clos\Range T^*_0 = \h\ominus \Ker T_0$  onto $\clos\Range T_0$.

We have  $\Ker\check T_0 =\Ker T_0\oplus\K$ and $(\Range\check T_0)^\perp =
(\Range T_0)^\perp\oplus\K$. Note that $\dim\Ker\check T_0 = \dim(\Range\check T_0)^\perp
= \infty$. Let $K:\Ker\check T_0 \longrightarrow\Ker\check T_0$ be a trace class nonnegative self-adjoint
contraction such that $\Ker K =\{\0\}$ and let $W$ be an isometry
from $\Ker\check T_0$ onto $(\Range\check T_0)^\perp$. Clearly, such an isometry  exists. Put
\begin{equation}
\label{2.12_oper_S}
Rf \df
\begin{cases}
T_0 f, & f \in \Ker(T_0)^\perp  = \h \ominus \Ker(T_0),  \\
WKf, & f \in \Ker(\check T_0) = \Ker(T_0)\oplus\K.
\end{cases}
\end{equation}

Let $f\in\h\oplus\K$ and $f = f_1 + f_2$, where   $f_1\in(\Ker T_0)^{\perp}$ and
$f_2\in\Ker\check T_0$. Then
   \begin{eqnarray*}
\|R f\|^2=\|T_0 f_1+W K f_2\|^2  = \|T_0 f_1\|^2 + 2 \RE(T_0 f_1, W K f_2) + \|WK f_2\|^2  \\
= \|T_0 f_1\|^2 +  \|WK f_2\|^2  \le\|f_1\|^2+\|f_2\|^2 = \|f\|^2.
   \end{eqnarray*}
Thus, $R$ is a contraction on $\h\oplus\K$. Moreover, $\Ker R=\{\0\}$ because $\Ker(T_0)=\Ker(W
K)=\{\0\}.$
It is easy to see that  $R$ has dense range in $\h\oplus\K$. Hence, $0 \in \rho(R)
\cup \gs_{\rm c}(R)$.  Let $ P_{\Ker\check T_0}$   be  the orthogonal projection onto
the subspace $\Ker\check T_0$. Combining \eqref{2.10} with \eqref{2.12_oper_S} and keeping in mind
that $K\in\bS_1$, we find that
$
R -\check T_0 = WK P_{\Ker\check T'_0}  \in  \bS_1.
$
It follows from what we have established in (ii) that there exists a spectral shift function $\bs\xi_1$  for the pair $\{\check T_0, R\}$, i.e., the identity
\begin{equation}
\label{2.14-trace_T-0,S}
\trace\big((\check T_0 - \l I)^{-1} - (R - \l I)^{-1}\big) = 
-\int_\T
\bs\xi_1(\z)(\z - \l)^{-2}d\z, \quad|\l|>1,
\end{equation}
holds.  Similarly,
 \begin{equation}
 \label{2.12_dif-ce_T'-S}
(\check T_1 - R)f =
\begin{cases}
(T- T_0) f, & f \in (\Ker T_0)^\perp,    \\
- WKf, & f \in \Ker\check T_0.
\end{cases}
\end{equation}
Hence, $\check T_1- R\in \bS_1$ and  there exists a spectral shift function
$\bs\xi_2$ for $\{\check T_1, R\}$
such that the following trace formula holds
\begin{equation}
\label{2.16-trace_T',S}
\trace\big((\check T_1 -\l I)^{-1} - (R - \l I)^{-1}\big) = 
-\int_\T
\bs\xi_2(\z)(\z - \l)^{-2}d\z, 
\quad |\l|>1 .
\end{equation}
Combining \eqref{2.14-trace_T-0,S} with \eqref{2.16-trace_T',S} and setting
$\bs\xi= \bs\xi_2 - \bs\xi_1$,  we obtain
\begin{equation}
\label{eq:2.18}
\trace\big((\check T_1-\l I)^{-1} - (\check T_0-\l I)^{-1}\big) = 
 -\int_\T \bs\xi(\z)(\z-\l)^{-2}d\z, \quad
|\l|>1.
\end{equation}
It follows from  \eqref{2.10} that
\begin{align*}
(\check T_1 - \l I)^{-1} - (\check T_0 -\l I)^{-1}& =  (T_1-\l I)^{-1}\oplus
-{\l}^{-1}I - (T_0-\l I)^{-1}
\oplus -{\l}^{-1}I \\
&=  \left((T_1-\l I)^{-1} - (T_0-\l I)^{-1}\right) \oplus \0.
\end{align*}
Finally, combining this identity  with  \eqref{eq:2.18}, we get trace formula
\eqref{eq:2.2}.

Suppose now that $0\in  \s_{\rm p}(T_0^*)\setminus  \s_{\rm p}(T_0)$. The above reasoning shows
that identity \eqref{eq:2.18}   holds with $\{\check T_1,\check T_0\}$  and $\bs\xi$ replaced
by $\{\check T_1^*,\check T_0^*\}$ and a function $\bs\xi_*$ in 
$\mathscr S\{\check T_1^*,\check T_0^*\}$.  Setting
$\bs\xi\df\overline{\bs\xi}_*$ and taking to the complex conjugate in this identity, we
arrive at \eqref{eq:2.18}.

(iv). Suppose now that $T_0$ is unitary. Then $0\in\rho(T_0)$ and in accordance with (ii),
 there exists $\vk\in\T$ such that $-\ov{\vk}\notin \sigma_p(T_0)\cup
\sigma_p(T)$. Therefore $L_0=\ri(I-\vk T_0)(I+\vk T_0)^{-1}$
and $L_1=\ri(I - \vk T_1) (I+\vk T)^{-1}$ are well defined.   Moreover, $L_0 =
L^*_0$,  $L_1$ is m-dissipative, and $\ri\in\rho(H_0)$. By  Theorem
\ref{thm:2.1} there exists a function $\bs\o$ in
$\mathscr S\{L_1, L_0\}$ such that $\im\bs\o \ge 0$. On the other hand,  in arcordance with (i), there is a spectral shift function
$\widehat{\bs\xi}$ such that trace formula \eqref{2.8_trace} holds, where
$\widehat{T}= \varkappa T$ and $\widehat{T}_0= \vk T_0$. By \eqref{2.6_connection_of_two_SSF},
 we get $\im\widehat{\bs\xi}\ge 0$ on $\bT$. Setting $\bs\xi(\z)=\widehat{\bs\xi}(\vk\z)$, $\z\in\T$,
we find that $\bs\xi$ is a spectral shift function for $\{T_1,T_0\}$ and $\im\bs\xi\ge 0$ on $\bT$. The case when $T_1$ unitary is treated
similarly. $\bl$

\medskip

Another approach to trace formulae for the resolvents for a pair $\{T_1, T_0\}$ with unitary
$T_0$ and a contraction $T_1$ satisfying  $T_1- T_0 \in \bS_1$ and $D_{T_1} \df (I-
T_1^*T_1)^{1/2}\in \bS_1$ was proposed by  A. Rybkin \cite{Ryb87,Ryb89,Ryb94}. His
formulae involve a complex spectral shift function, which is (A)-integrable but not Lebesgue integrable.

\section{Real-valued spectral shift functions}

\medskip

We denote by $\mathscr S\{T_1,T_0\}$ the set of  all spectral shift functions for the pair $\{T_0,T_1\}$.
Let $\bs\xi\in \mathscr S\{T_1,T_0\}$. Clearly, an $L^1$ function  $\bs\xi_\natural$
belongs to $\mathscr S\{T_1,T_0\}$ if and only if

  \begin{equation}
  \label{8.0-dif-ce_of_two_SSF}
\bs\xi_\natural-\bs\xi\in H^1.
   \end{equation}
   
  \begin{thm}
  \label{Zigmund_cond}
Let $\{T_1,T_0\}$  be a pair of contractions satisfying $T_1-T_0 \in
\gotS_1$ and let $\bs\xi\in \mathscr S\{T_1,T_0\}$. 
Suppose that
\begin{equation}
\label{2.36}
(\im\bs\xi)\log(1+|\im\bs\xi|)\in L^1(\T).
    \end{equation}
Then there exists a real function in $\mathscr S\{T_1,T_0\}$.

If   $\im\bs\xi \ge 0$ for a.e. on $\T$ and there is a real function in
$\mathscr S\{T_1,T_0\}$, then \eqref{2.36}
holds.
\end{thm}

\Pf For $g\in L^1(\T)$, we put $\pp_+g\df\sum_{j\ge0}\widehat g(j)z^j$ and 
$\pp_-g\df g-\pp_+g$. 
We have
$$
\bs\xi=\re\bs\xi+{\rm i}\im\bs\xi
=\re\bs\xi+{\rm i}\pp_-(\im\bs\xi)+{\rm i}\pp_+(\im\bs\xi).
$$
By Zygmund's theorem (see \cite{Koo80}, Sect. V.C.3), inequality \rf{2.36} implies that both $\pp_+(\im\bs\xi)$ and $\pp_-(\im\bs\xi)$ belong to $L^1(\T)$. Put
\bay
\label{ksifl}
\bs\xi_\flat\df\re\bs\xi+{\rm i}\pp_-(\im\bs\xi)-{\rm i}\ov{\pp_-(\im\bs\xi)}.
\ey
Clearly, $\bs\xi_\flat$ is a real function in $L^1(\T)$ and
$$
\bs\xi-\bs\xi_\flat={\rm i}\big(\pp_+(\im\bs\xi)+\ov{\pp_-(\im\xi)}\,\big)\in H^1.
$$
By \rf{8.0-dif-ce_of_two_SSF}, $\bs\xi_\flat\in\mathscr S\{T_1,T_0\}$.

Suppose now that $\bs\xi$ and $\bs\xi_\sharp$ belong to $\mathscr S\{T_1,T_0\}$,
$\im\bs\xi\ge0$  and $\bs\xi_\sharp$ is real. Then by \rf{8.0-dif-ce_of_two_SSF}, $\bs\xi-\bs\xi_\sharp\in H^1$ and
$\im(\bs\xi-\bs\xi_\sharp)=\im\bs\xi\ge0$. It follows from M. Riesz' theorem (see 
\cite{Koo80}, Sect. V.C.4) that \rf{2.36} holds. $\bl$

   \begin{thm}
   \label{prop_real_SSF_Muckenh}
Let $T_0$ and $T_1$ be contrarctions with trace class difference. 
Suppose that $\bs\xi\in\mathscr S\{T_1,T_0\}$ and $\im\bs\xi\in L^p(\T,w)$
for a weight $w$ satisfying the Muckenhoupt condition $({\rm A})_p$, $1<p<\be$.
Then there is a real function 
$\bs\xi_\flat\in\mathscr S\{T_1,T_0\}\cap L^p(\T,w)$.
\end{thm}

\Pf It follows easily from H\"older's inequality that $L^p(\T,w)\subset L^1(\T)$.
We can take the function $\bs\xi_\flat$ defined by \rf{ksifl}, 
repeat the argument in the proof of Theorem \ref{Zigmund_cond}
and use the Hunt--Muckenhoupt--Wheeden theorem (see \cite{Koo80}) instead of Zygmund's theorem. $\bl$

\medskip
  
Consider now the special case of pairs $\{L_1,L_0\}$, where $L_0 =\0$, $L_1= L_0+{\rm i}V$ and $\0\le V \in\bS_1$. Clealy, $L_1$  is bounded, and so it is m-dissipative.

\begin{thm}
\label{criterion_for_exist_real_SSF}
Let $L_0 = \0$,  $L_1 = L_0 + {\rm i}V$,   $\0 \le V \in \bS_1$, and let $\{\ga_n\}_{n\ge0}$ be the sequence of eigenvalues of $V$
counted with multiplicities. 
Then
\begin{itemize}
\item [(i)]
The set  $\mathscr S\{L_1,L_0\}$  contains  a purely imaginary  function;

\item  [(ii)]
 The  set  $\mathscr S\{L_1,L_0\}$  contains  a real-valued  function if and only if  
  \begin{equation}
  \label{Zigmund_cond_for_eigen}
\sum_{n\ge0}\ga_n |\log\ga_n| < \infty.
\end{equation}
\end{itemize}
  \end{thm}

\Pf
(i)  Consider the perturbation determinant $\gD_{L_1/L_0}$ defined  by
\begin{equation}
\label{3.12_pert_deter}
\gD_{L_1/L_0}(z) = \det(I +\ri V(L_0-zI)^{-1}) = \det(I - \ri z^{-1}V), \quad z \in \bC
\setminus \{0\}.
\end{equation}
The operator $V = V^*$ being of trace class, admits a spectral decomposition
\begin{equation}\label{3.13_spec_decom-n}
V = \sum_{n\ge0}\ga_nP_n = \sum_{n\ge0}\ga_n(\cdot,e_n)e_n, \qquad \ga_n \ge1,\quad n\in \bN,
\end{equation}
where  $\{e_n\}^\infty_{n=1}$ is an orthonormal basis, the orthogonal projection 
$P_n$ is defined by
$P_n=(\cdot,e_n)e_n$ is   , and $C_0\df
\sum^\infty_{n=1}\ga_n < \infty.$
Setting $V_n= \ga_nP_n = \ga_n(\cdot,e_n)e_n$, $L_n = L_0 - \ri V_n = -\ri V_n$  and
combining \eqref{3.12_pert_deter} with \eqref{3.13_spec_decom-n}, we obtain
   \begin{equation}
   \label{3.13-determ-t}
\gD_{L_1/L_0}(z) = \prod_{n\ge0}(1 - \ri\ga_n z^{-1}) =
\lim_{N\to\infty}\prod^N_{n=0}(1 -\ri\ga_n z^{-1}).
   \end{equation}
By Lemma 5.5, (i) of \cite{MalNei2015}, there is a non-negative  function $\eta_n$ in
$L^1(\bR)$ such that
\begin{equation}
\label{3.14-u-n}
\begin{split}
u_n(z) &
\df\gD_{L_n/L_{n-1}}(z)  \df \det(I + \ri V_n (L_{n-1}-zI)^{-1}) \\
& = 1 - \ri\ga_n z^{-1} =\exp\Big(\ri\pi^{-1}\int_\bR \eta_n(t)(t-z)^{-1}dt\Big),
\qquad z\in \bC_-,
\end{split}
    \end{equation}
where
\begin{equation}
\label{3.15-eta_n}
\eta_n(t) = \lim_{y\uparrow 0} \log(|u_n(t - iy)|)  = \log(|1 - \ri\ga_n t^{-1}|)  =
 (1/2)\log\left(1 +\ga^2_nt^{-2}\right) \ge 0.
\end{equation}
Integrating by parts,  we see that $\eta_n\in L^1(\bR)$ and
$$
\sum_{n\ge0} \int_\bR \eta_n(t) = 
\sum_{n\ge0} \int^\infty_0\log\left(1 + \tfrac{\ga^2_n}{t^2}\right)dt =
2\sum_{n\ge0} \int^\infty_0\frac {\ga^2_n}{t^2 + \ga^2_n}\, dt =  \pi
\sum_{n\ge0} \ga_n = \pi C_0.
$$
By the B. Levi theorem,   the function  $\eta$,
   \begin{equation}
   \label{3.17-density-eta-n}
\eta(t) \df \sum_{n\ge0}\eta_n(t) =(1/2)\sum_{n\ge0}\log\left(1 +\ga^2_nt^{-2}\right),
  \end{equation}
is well defined for almost all $t\in \bR$, non-negative, and  $\eta\in
L^1(\bR)$.  Combining  relations \eqref{3.13-determ-t}--\eqref{3.17-density-eta-n}, we
arrive at the representation
\begin{equation}
\gD_{L_1/L_0}(z) = \exp\Big(\ri\pi^{-1}\int_\bR \eta(x)(x-z)^{-1}dx\Big) \qquad
z\in \bC_-.
\end{equation}
Taking the logarithmic derivative one arrives at trace formula \eqref{eq:2.1} showing
that $\ri\eta$ is a (purely imaginary)  spectral shift functions for the pair $\{L_1,L_0\}$, i.e.,
$\ri\eta\in \mathscr S\{L_1,L_0\}$.

(ii)
It follows from \eqref{3.15-eta_n} that the harmonic conjugate to ${\eta}_n$
function $\widetilde{\eta}_n$  is given by
$\widetilde{\eta}_n(t)= \arctan(\a_n/t),\  t \in \bR, \  n\ge0$,
and
$$
\widetilde\eta(t)  = -\lim_{y\uparrow0}\frac{1}{\pi}\int_\bR\frac{\eta(x)(x-t)}{(x-t)^2
+ y^2}dx = \sum_{j\ge0}\widetilde\eta_j(t) =
\sum_{j\ge0}\arctan\frac{\alpha_n}{t},\quad t \in \bR,
$$
where the convergence is pointwise. Clearly, $\widetilde\eta$ is real-valued and
$\widetilde\eta+\ri\eta\in H^1(\bR, (1+t^2)^{-1})$ if and only if
$\widetilde\eta\in L^1(\R,(1+t^2)^{-1})$. Therefore, by Proposition
3.8 of \cite{MalNei2015},  a real spectral shift function for $\{L_1,L_0\}$ exists if and only if $\widetilde
\eta \in L^1(\R,(1+t^2)^{-1})$. Thus, it remains to show that the latter
inclusion is equivalent to condition \eqref{Zigmund_cond_for_eigen}.

Making use the change of variables $x = \ga_j/t$ we obtain
  \begin{equation}\label{eq:3.48}
\begin{split}
& \int_\bR \frac{|\widetilde\eta(t)|}{1+t^2}\,dt = \sum_{j=1}^\infty \int_\bR
\frac{|\widetilde\eta_j(t)|}{1+t^2}\,dt = 2\sum_{j=1}^\infty \int^\infty_0
\frac{\widetilde\eta_j(t)}{1+t^2}dt = 2\sum_{j=1}^\infty\int^\infty_0
\frac{\arctan(\ga_j/t)}{1+t^2}\,dt  \\
& = 2\sum_{j=1}^\infty \ga_j\int^\infty_0\frac{\arctan x}{\ga^2_j + x^2}dx = 2
\sum_{j=1}^\infty \ga_j\int^1_0\frac{\arctan x}{\ga^2_j+x^2}\,dx + 2\sum_{j=1}^\infty
\ga_j\int^\infty_1\frac{\arctan x }{\ga^2_j+x^2}\,dx.
\end{split}
   \end{equation}
It is easily seen that for some $c_0>0$,
  \begin{equation*}
c_0 x\le\arctan x\le x, \quad x\in[0,1], \quad \text{and}\quad
 \int^{\infty}_1(\arctan x)(\alpha_j^2 + x^2)^{-1}dx \le \frac{\pi}{2}.
 \end{equation*}
Combining these estimates with  \eqref{eq:3.48} leads to two-sided estimate
   \begin{equation}
c_0 \sum_j\alpha_j\int^1_0\frac{2x}{x^2+\alpha^2_j}\,dx  \le \int_\bR
\frac{|\widetilde\eta(t)|}{1+t^2}\,dt  
 \le  \sum_j\alpha_j\int^1_0\frac{2x}{x^2+\alpha^2_j}\,dx + \pi C_1,
  \end{equation}
where $C_1 = \sum_{j\ge0}\ga_j$. To complete the proof, it remains to observe that
$$
\int^1_0 2x(\ga^2_j+x^2)^{-1}dx  = \log\left((1 + \ga^2_j)\ga_j^{-2}\right) \sim
-2\log\alpha_j = 2|\log\alpha_j|\qquad \text{as} \quad  j\to\infty.\quad\bl
$$

Using the Cayley transform, we deduce the following result.

  \begin{cor}
  \label{cor_zigmund}
Let $T_0=I$ and let $T_1$ be a  self-adjoint contraction such that $I-T_1\in\ \bS_1$ and
let $\{\lambda_n\}_{n\ge0}$ be the sequence of eigenvalues of $T_1$ counted with multiplicities. If\  $\Ker(I \pm T) =
\{\0\}$, then $\mathscr S\{T_1,T_0\}$ contains no real function if and only if
  \begin{equation}
  \label{Zigmund_cond_for_contr-s}
\sum_n(1-\lambda_n)(1+\lambda_n)^{-1}\big|\log\big((1-\lambda_n)(1+\lambda_n)^{-1}\big)\big|
 = \infty.
\end{equation}
\end{cor}

\Pf
Let $L_0$ and $L_1$ be the  Cayley transforms of $T_0$ and $T_1$ (see
\eqref{2.4_oper_H_and_H_0}). Clearly,  $L_0=\0$ and $L_1=-\ri V$, where
$V=(I-T_1)(I+T_1)^{-1}\in\bS_1$ and $V\ge 0$. Therefore condition
\eqref{Zigmund_cond_for_eigen} is transformed into condition
 \eqref{Zigmund_cond_for_contr-s}. It
remains to apply Theorem \ref{criterion_for_exist_real_SSF} and take formula
\eqref{2.6_connection_of_two_SSF}  into account. $\bl$

\medskip

{\bf Remark.}
In connection with Theorem \ref{main_theor} and  Corollary \ref{cor_zigmund}  we mention
Theorem 3.3  from \cite{AN}. It says that the class $\mathscr S\{T_1,T_0\}$ contains
a real function whenever in addition to the condition $T_1-T_0 \in \bS_1$ the
following conditions are satisfied: $I - |T_1| \in \bS^0_1$ and $I - |T^*_0| \in
\bS^0_1$, where
$
\bS^0_1 \df\big\{A :~ \sum_j
s_j(A)\big|\log s_j(A)\big| < \infty\big\} \subset \bS_1
$
and $\{s_j(A)\}_{j\ge0}$  is the sequence of singular values of an  operator $A$. Theorem \ref{criterion_for_exist_real_SSF}  shows that Theorem
3.3 of \cite{AN} is sharp.

In \cite{MakScZ1} the authors studied spectral shift functions for a  pair $\{L_0,L_0-\ri V\}$ with $L_0=L^*_0$ and  $0\le V\in
\bS_1$. In particular,  it was shown that $w_2(t) = \lim_{y\downarrow
0}\arg\bigl(\det_{L_1/L_0}(t+\ri y)\bigr)$ is in $L^1(\Bbb R)$ if and only if the operator $L_0-\ri V$ is self-adjoint.

\

\section{\bf The main results}

\medskip

In this section we establish the main results of the paper. The following theorem 

\begin{thm}
\label{szhatiya}
Let $T_1$ and $T_0$ be contractions such that $T_1-T_0\in\bS_1$ and let $f\in\OLA$.
Then $f(T_1)-f(T_0)\in\bS_1$ and the following trace formula holds
\bay
\label{okfsl}
\trace\big(f(T_1)-f(T_0)\big)=\int_\T f'(\z)\bs\xi(\z)\,d\z
\ey
for an arbitrary spectral shift function $\bs\xi$ for the pair $\{T_1,T_0\}$.
In particular, 
\bay
\label{intksi}
\int_\T\bs\xi(\z)\,d\z=\trace(T_1-T_0).
\ey
Conversely, if for a function $f$ in $\CA$, $f(T_1)-f(T_0)\in\bS_1$ whenever $T_1$ and $T_0$ are contractions such that $T_1-T_0\in\bS_1$, then $f\in\OLA$.
\end{thm}

\Pf By Theorem \ref{fsdm}, formula \rf{formsledlyaszha} holds for some complex Borel measure $\nu$. On the other hands, by Theorem \ref{main_theor}, if 
$\bs\xi\in\mathscr S\{T_1,T_0\}$, then formula \rf{okfsl} in the case 
$f(z)=(z-\l)^{-1}$, $|\l|>1$. Since the linear combinations of such functions are dense in the class $X\df\{f:~f'\in H^1\}$, it follows that $\int_\T f'(\z)\bs\xi(\z)\,d\z=\int_\T f'(\z)\,d\nu(\z)$ for every $f$ in $X$. This implies the following equality for the Fourier coefficients: $\widehat{\bs\xi}(j)=\widehat\nu(j)$, $j<0$. By the brothers Riesz theorem (see \cite{Koo80}, Ch. II, Sect. A), $\nu$ is absolutely continuous with respect to Lebesgue measure and its Radon-Nikodym density is a spectral shift functions for the pair $\{T_1,T_0\}$. This implies formula \rf{okfsl} for an arbitrary function $f$ in $\OLA$. To get formula \rf{intksi}, it suffices to put $f(z)=z$.

The last conclusion of the theorem is an immediate consequence of the corresponding fact for functions of unitary operators, see \cite{AP+} and the fact that 
$\OLA=\OL_\T\bigcap C_{\rm A}$ (see \cite{KS} and \cite{AP}). $\bl$

\begin{thm}
\label{dissip}
Let $\{L_1,L_0\}$ be a pair of resolvent comparable maximal dissipative operators
and let $f$ be a function analytic in $\C_+$ and such that the function
$$
\z\mapsto f\big(({\rm i}-\z)({\rm i}+\z)^{-1}\big),\quad\z\in\C_+,
$$
belongs to $\OLA$.
Then $f(L_1)-f(L_0)\in\bS_1$ and
$$
\trace\big(f(L_1)-f(L_0)\big)=\int_\R f'(t)\bs{\o}(t)\,dt,
$$
where $\bs\o$ is a spectral shift function for the pair $\{L_1,L_0\}$.
\end{thm}

\Pf The theorem can be proved by passing to the contractions $T_0=(\ri I-L_0)(\ri I+L_0)^{-1}$
and $T_1=(\ri I-L_1)(\ri I+L_1)^{-1}$ and applying Theorem \ref{szhatiya}. $\bl$

\medskip

Note that Theorem \ref{dissip} eliminates the additional assumption $(\rho(L_0) \cup
\sigma_{\rm c}(L_0))\cap \bC_+ \not= \varnothing$ in the statement of Theorem 
\ref{thm:2.1}.

\begin{appendix}

\section{\bf Pairs of  accumulative  and dissipative  operators}\label{sec.V.2}

\medskip

{\bf A1.~Auxiliary results.}
Here we improve Theorem 3.13 from \cite{MalNei2015} mentioned in the Introduction. To
this end we need several auxiliary results.

Recall (see \cite{SNF}) that a bounded holomorphic 
 function $W: \bC_+ \to \cB(\h)$ is called an {\it outer function}
if multiplication by $W$ on $H^2_{\h}(\bC_+)$ has 
dense range. Recall that the strong limit $W(t)\df
\lim_{y\downarrow 0}W(t +\ri y)$ exists for almost all $t\in\R$, see \cite{SNF}.

  \begin{lem}
  \label{lemA1_outer_func}
Let $W: \bC_+ \to \cB(\h)$ be a contractive holomorphic function.  If $\re W(z)>\0$ for $z\in \bC_+$ and $\Ker\re W(t)=\{\0\}$
 a.e. on $\bR $, then $W$ is an outer function.
  \end{lem}

\Pf
Assume the contrary. Then there exists $g\in H^2_{\h}(\C_+)$ such that
  \begin{equation}
  \label{Append_1}
(Wf,g)_{H^2_{\h}} = \int_{-\infty} ^{\infty} \bigl(W(t)f(t),
g(t)\bigr)_{\h}\,dt = 0,\quad \text{for all}\quad f\in H^2_{\h}(\bC_+).
    \end{equation}
Putting $f=g$ in \eqref{Append_1}, taking the real part, and observing that $\re\Theta(t)
 \ge 0$  for a.e. $t\in \bR$, yields
$
\int_\R \|(\re W(t))^{1/2}g(t)\|^2_{\h}\,dt = 0$, and so
$(\re W(t))g(t)=0$ a.e. on $\bR$.
Hence, $g(t)=0$  a.e. and multiplication by $W$ has dense range in
$H^2_{\h}(\bC_+)$. $\bl$

\medskip

{\bf Definition.} 
A holomorphic function $F:\bC_{+}\cup\bC_{-}\to \mathcal B(\h)$ is called  an
{\it $R$-function} (in short $F\in R[\h]$)  if $\IM z\cdot \IM F(z)>0$ 
and $F(z)^*=F(\overline z)$ for $z\in\bC_{\pm}$. We also use the notation
$R^{\rm u}[\h]$ for $\{F\in R[\cH]:0\in\rho\bigl(\im F(z)\bigr),\ z\in\bC_{\pm}\}$.

\begin{lem}
\label{IV.3a}
Let $M\in R^{\rm u}[\h]$  and $B$ a bounded accumulative operator on $\h$. Then
\item[\;\;\rm (i)] If $0 \le V_+ \le |\im B| = -\im B\ge 0$ and $V_+ \in \bS_1$, then  the function $w_+\lb\df \det(I +
\ri V_+(B - M)^{-1})$  is holomorphic and contractive in $\bC_+$. Moreover,
$w_+$ is an outer function, and so it admits a representation
  \begin{equation}
  \label{4.180}
w_+(z) = \vk_+ \exp\left(\frac\ri\pi
\int_\R ((t-z)^{-1} -t(1+t^2)^{-1}) 
\eta_+(t)dt\right), \qquad z \in \bC_+,
\end{equation}
where $\eta_+(t) = -\log|w_+(t)|\ge 0$  a.e. on $\bR$,  $\eta_+ \in
L^1(\R,(1+t^2)^{-1}),$ and   $\varkappa_+ \in \bT$.

\item[\;\;\rm (ii)] If $V \le |\im B| = -\im B$ and $V = V^* \in
\bS_1$, then $w\df\det(I + \ri V(B - M)^{-1})$  is an
\emph{outer function} in  $\bC_+$ and
 admits a representation
\begin{equation}
\label{4.190}
w(z) = \vk \exp\left(\frac{\ri}{\pi}
\int_\R ((t-z)^{-1} -t(1+t^2)^{-1})  
\eta(t)dt\right), \quad  z \in \C_+,\quad \varkappa_+ \in \bT,
\end{equation}
with a
 real function $\eta\in L^1(\R,(1+t^2)^{-1})$. Besides, $\eta(t) =
-\log|w(t)|$  a.e. on $\R.$
\end{lem}

\Pf
(i)  Put
\begin{equation}
\label{fun-n_W_+}
W_+(z)\df I_{\h} +\ri{V_+}^{1/2}(B - M(z))^{-1}{V_+}^{1/2}, \quad z \in \C_+.
\end{equation}
First, we show that  $W_+$ is a holomorphic contractive $\h$-valued function in
$\C_+$. Since $\im M(z)\ge \varepsilon(z)I$, $z \in \C_+$,  for a positive function $\e$ on $\C_+$ and $\im B\le\0$, the function
$(B - M)^{-1}$ is well defined and holomorphic in $\C_+$.  Besides, taking
into account the inequality $(\im B + V_+) \le\0$,  we get
$$
I- W_+^*(z)W_+(z)
= V_+^{1/2}(B^* - M^*(z))^{-1}(2\im M(z) - 2\im B - V_+)(B - M(z))^{-1}V_+^{1/2} \ge\0.
$$
Hence, $W_+$  is contractive. Next, we show that $W_+$ is  an outer
$\mB(\h)$-valued function in $\C_+$. To this end we set $C\df B + \ri V_+$ and
  \begin{equation}
  \label{fun-n_W_-}
W_-(z) \df W_+(z)^{-1} = I_{\cH} -\ri V_+^{1/2}(C - M(z))^{-1}V_+^{1/2}\ , \quad z \in
\C_+.
  \end{equation}
Besides, keeping in mind that $\im C =\im B + V_+ \le 0$, one can easily obtain
  \begin{equation}
  \label{Real-part_of_W_-}
\RE W_-(z)= I_{\h} + V_+^{1/2}(C - M(z))^{-1}(-\im C+\im M(z)(C^*-M^*(z))^{-1}V_+^{1/2}\ge I_{\cH} 
  \end{equation}
for  $z \in \C_+$.   It follows with account of \eqref{fun-n_W_-} that
  \begin{equation}
  \label{Real-part_of_W_+}
\RE W_+(z) =   W_-(z)^{-1}(\RE W_-(z))(W_-(z)^{-1})^* \ge W_-(z)^{-1}(W_-(z)^{-1})^* >\0.
  \end{equation}
Since $V \in \bS_1$, the limits $W_{\pm}(x)=\lim_{y\downarrow 0}W_{\pm}(x+\ri y)$ exist
in the $\bS_2$-norm  for a.e. $x\in\bR$. Therefore passing to the limit in the identity
$W_{\pm}(x+\ri y)W_{\mp}(x+\ri y)=I$ as $y\downarrow 0$, we obtain $W_{\pm}(x)W_{\mp}(x)=I$ 
a.e. on $\bR$.  Passing to the limit as $y\downarrow 0$ in \eqref{Real-part_of_W_+} we obtain
  \begin{equation*}
  \label{Real-part_of_W_+(x)}
\RE W_+(x) = W_-(x)^{-1}(\RE W_-(x))(W_-(x)^{-1})^* \ge W_-(x)^{-1}(W_-(x)^{-1})^* >
\0\quad\mbox{a.e.}.
  \end{equation*}
Hence, $\Ker \RE W_+(x)= \{0\}$  a.e. on $\bR$.  By Lemma \ref{lemA1_outer_func},
$W_+$ is an outer function in $\bC_+$ and so is $w_+(z)= \det W_+(z)$ (see  \cite{SNF}). Being contractive and outer in $\C_+$, the function $w_+$ admits
a representation \eqref{4.180}  (see \cite{Koo80}, Sect. VI C).

(ii) Let $V = V_+ - V_-$, $V_\pm \ge 0$. Put $B_-\df B -\ri V_-$. Since 
$\im B_- =\im B-V_-\le\0$, the operator $B_-$ is accumulative. One can easily check that for $z \in \C_+$,
\begin{displaymath}
\det(I +\ri V(B - M(z))^{-1}) =
{\det(I +\ri V_+(B_- - M(z))^{-1})}(\det(I +\ri V_-(B_- -
M(z))^{-1}))^{-1}.
\end{displaymath}
The assumption $V \le -\im B$ yields  $\0 \le V_+ \le -\im B + V_- = -\im B_-$. In accordance
with 
(i) there exist a non-negative $\eta_+ \in L^1(\R,(1+t^2)^{-1})$ and
$\varkappa_+ \in \T$  such that
  \begin{displaymath}
\det(I +\ri V_+(B_- - M(z))^{-1}) = \vk_+ \exp\left(\frac\ri\pi
\int_\R ((t-z)^{-1} -t(1+t^2)^{-1}) 
\eta_+(t)dt\right),\quad z \in \C_+.
\end{displaymath}
Similarly, since  $\0 \le V_- \le -\im B + V_- = -\im B_-$ and (i) ensures a
representation
\begin{displaymath}
\det(I +\ri V_-(B_- - M(z))^{-1}) = \vk_- \exp\left(\frac\ri\pi
\int_\R ((t-z)^{-1} -t(1+t^2)^{-1}) 
\eta_-(t)dt\right),\ z \in \C_+
\end{displaymath}
for a nonnegative function $\eta_- \in L^1(\R,(1+t^2)^{-1})$ and $\vk_-
\in \T$.
Setting $\vk \df\vk_+/\vk_- \in \T$ and $\eta\df\eta_+-\eta_-$,  we arrive at  \eqref{4.190}. $\bl$

\medskip

The following lemma generalizes  Krein's result \cite{Kr} on exponential representation of the
perturbation determinant of a pair of selfadjoint operators.

\begin{lem}
\label{IV.2}
  Let $M\in R^{\rm u}[\h]$,  let $B$ a bounded accumulative (possibly,
self-adjoint) operator on $\h$
and let  $V = V^* \in \bS_1$.  Then there exist  a \emph{real-valued function}
$\xi$ in $L^1(\R,(1+t^2)^{-1})$ and a constant $c >0$   such that the  equality
   \begin{equation}
   \label{4.4a}
\det\left(I + V(B - M(z))^{-1}\right)= c\, 
\exp\left(\frac1\pi
\int_\R ((t-z)^{-1} -t(1+t^2)^{-1}) 
\xi(t)dt\right)
  \end{equation}
holds for $z \in \C_+$. Moreover, if $V = V_+ \ge 0$, then $\xi= \xi_+\ge 0$.
  \end{lem}

\Pf
(i) Assume first that $V = V_+ \ge\0$.  Consider the operator-valued Nevanlinna function $\O_+$, 
$
\gO_+(z) \df I + {V_+}^{1/2}(B - M(z))^{-1}{V_+}^{1/2}, ~ z \in \C_+.
$ 
Since $\gO_+(z)$ is $m$-dissipative for $z \in \C_+$ and $0 \in \rho(\gO_+(z))$, $z \in
\C_+$, the operator-valued function $\log(\gO_+(z))$ is well-defined 
 for $z \in \C_+$ (see \cite{GMN99}).  Besides, since $V_+ \in \bS_1$, we have
$\log(\gO_+(z)) \in \bS_1$ and Theorem 2.8 of \cite{GMN99}) ensures the existence of
a non-negative measurable function $\Xi_+:
\R \to \bS_1$ such that the following representation holds
  \begin{equation*}
\log(\gO_+(z)) = \gO_+ + \frac{1}{\pi}\int_\R \left(\frac{1}{t-z} -
\frac{t}{1+t^2}\right)\Xi_+(t)dt, \quad z \in \C_+,\quad \int_\R
\frac{\Xi_+(t)}{1+t^2}\,dt\in \cB(\h).
  \end{equation*}
Here the integral is understood   in the  weak sense and $\gO_+ = \gO_+^* \in \bS_1$.
Taking the traces in this identity and setting $\xi_+(t)\df \trace(\Xi_+(t))$, $t \in \R$, we
arrive at the representation
  \begin{equation*}
\trace(\log(\gO_+(z))) = \trace(\gO_+) + \frac{1}{\pi}
\int_\R((t-z)^{-1} -
t(1+t^2)^{-1})\xi_+(t) dt, \quad z \in \C_+.
  \end{equation*}
Taking into account the identity $\det(I + G) = e^{\trace(\log(I+G))}$  valid for any
dissipative $G \in \gotS_1(\gotH)$,  we  arrive at  \eqref{4.4a} with $V=V_+$,
$\xi= \xi_{+},$ and  $c= c_+= \exp\{\tr(\gO_+)\} > 0$.

(ii) Passing to the general case, we start with  the spectral decomposition $V = V_+ -
V_-$ $(V_\pm \ge\0$, $V_{\pm}V_{\mp} =\0)$, and  set $B_-\df B - V_-$. It follows from
the identity
$$
\left(I + V(B - M(z))^{-1}\right) \left(I + V_-(B_- - M(z))^{-1}\right) =  I + V_+(B_- -
M(z))^{-1}
$$
that for
  \begin{equation}
  \label{4.20}
\det\left(I + V(B - M(z))^{-1}\right) = \frac{\det\left(I + V_+(B_- - M(z))^{-1}\right)}
{\det\left(I + V_-(B_- - M(z))^{-1}\right)},\qquad z \in \C_\pm.
  \end{equation}

Note that representations \eqref{4.4a} with $V_{\pm}\ge\0$ and $\xi_{\pm}\ge 0$
instead of $V$ and $\xi$  have  just  been established in (i).
Combining these representations with \eqref{4.20} and setting $\xi\df
\xi_{+}- \xi_{-}$, we arrive at \eqref{4.4a}. $\bl$

\medskip

{\bf A2.~Boundary triplets and Weyl functions.}
Let us  recall basic facts of the theory of boundary triplets. Let $A$ be a densely
defined symmetric operator in $\h$.

\medskip  

{\bf Definition.}
A triplet $\Pi = \{\K, \Gamma_0, \Gamma_1\}$, where $\K$ is an auxiliary Hilbert space
and $\Gamma_0,\Gamma_1:\ \dom(A^*)\rightarrow \K$ are linear operators,  is called a
{\it boundary triplet for} $A^*$ if "Green's identity"
    \begin{equation}
    \label{2.0}
(A^*f,g) - (f,A^*g) = (\gG_1f,\gG_0g)_{\K} - (\gG_0f,\gG_1g)_{\K}, \qquad
f,g\in\dom(A^*),
   \end{equation}
holds and the mapping $\gG=(\Gamma_0,\Gamma_1):  \dom(A^*) \rightarrow \cH \oplus
\cH$ is surjective.

\medskip

A boundary triplet $\Pi=\{\K,\gG_0,\gG_1\}$ for $A^*$ exists whenever $A$ has equal deficiency indices  $n_+(A) =
n_-(A)$. Note also that  $n_\pm(A) = \dim(\K)$.  With any boundary triplet $\Pi$ we can associate a self-adjoint extension $A_0$ of $A$,
$A_0\df A^*\big|\ker(\gG_0) = A_0^*$.

A closed extension $\wt A$ of $A$ is called {\it proper} if $A\subseteq \wt A \subseteq A^*.$ The
set of proper extensions is denoted by $\Ext_A$. A dissipative (accumulative)
extension of $A$ is always proper.

The role of a boundary triplet in extension theory is similar to the role of a coordinate system in
analytic geometry. In particular, it allows one to  parameterize the set $\Ext_A$ by
means of the set $\widetilde\cC(\cH)$  of closed  linear relations in $\K$ (i.e., subspaces
in $\K\oplus\K)$.

\begin{lem}[{\cite{DM1}}]
\label{prop2.1}
Let  $\Pi=\{\K,\gG_0,\gG_1\}$  be a boundary triplet for  $A^*.$ The map
\begin{equation}
\label{bij}
\Ext_A \ni \widetilde A \to\Theta\df\Gamma \dom(\widetilde A) =\{\{\Gamma_0 f,\Gamma_1f \} : \
f\in \dom(\widetilde A) \} \in \widetilde\cC(\K)
\end{equation}
establishes  a bijective correspondence between 
$\Ext_A$ and  the set of all linear relations in $\K$. We write  $\wt A = A_\Theta$ if
$\wt A$ corresponds to $\Theta$ via \eqref{bij}. Moreover, the following holds:

\item[\;\;\rm (i)]  $A_\Theta$ is symmetric (self-adjoint) if
  and only if  $\Theta$ is symmetric (self-adjoint).

\item[\;\;\rm (ii)]
  $A_\Theta$ is m-dissipative  (m-accumulative) if
  and only if so is $\Theta$.

\item[\;\;\rm (iii)] The extensions $A_\Theta$ and $A_0$ are
disjoint, i.e., $\dom A_W\cap \dom A_0 = \dom A$, 
if and only if $\Theta= {\rm graph}(B)$ and $B$ is a closed operator in $\K$,
and \eqref{bij} takes the form
\begin{equation} 
\label{2.2}
\wt A =A_B \df A_{{\rm graph}(B)} = A^*\big|\ker(\gG_1- B\gG_0).
\end{equation}
\end{lem}

The operator $B$ is called \emph{the boundary operator} of $\wt A$ with respect to the
triplet $\Pi$.

\medskip

{\bf Definition} (\cite{MN1}).
A boundary triplet $\Pi =\{\K,\gG_0,\gG_1\}$ for $A^*$ is called {\it regular} for the
family  $\{\wt A_j\}^N_{j=1}$ of proper extensions if there exist bounded  operators
$B_j \in \cB(\K)$  such that $\wt A_j = A_{B_j}\df A^*\big|\Ker(\gG_1-
B_j\gG_0),$\ $j\in \{1, \ldots, N\}$, $($cf. $\eqref{2.2})$.

\begin{thm} 
{\em(\cite{MN1}, Th. 3.5).}
\label{prop3.3}
Let $A$ be as above  and let $\{\wt A_j\}^N_{j=1} \subset \Ext_A$ and 
$\wt A_{N+1}=\wt A \in \Ext_A$.  Assume also that $\bigcap^{N+1}_{j=1}\rho(\wt A_j) \not=
\varnothing$ and
  \begin{equation}
  \label{7.7A}
(\wt A - z_1 I)^{-1} - (\wt A_j - z_1 I)^{-1} \in \bS_\infty,\quad z_1 \in
\bigcap^{N+1}_{k=1}\rho(\wt A_k), \quad j \in \{1,\ldots,N\}.
  \end{equation}
If there is  $z_2 \in \rho(\wt A)\cup\s_{\rm c}(\wt A)$ such that $\IM(z_1)\IM(z_2) < 0$,
then the family  $\{\wt A_j\}^{N+1}_{j=1}$   admits a regular boundary triplet. In
particular, the last condition is satisfied if  $\rho(\wt A)\cap \C_{\pm}\not
=\varnothing$.
\end{thm}

The main analytical tool in this theory is played by the so-called  Weyl function.  

\medskip

{\bf Definition} ({\cite{DM1}}).
Let $\Pi=\{\K,\gG_0,\gG_1\}$ be a boundary triplet  for $A^*$,
$A_0=A^*\big|\ker(\gG_0)$, and $\mathfrak N_z \df\ker(A^* - z I)$. The
operator valued functions $\gamma: \rho(A_0)\rightarrow  \cB(\K,\h)$ and
$M :\rho(A_0)\rightarrow \cB(\cH)$ defined by
\begin{equation}
\label{2.3A}
\gamma(z)=\bigl(\Gamma_0\big|\mathfrak N_z\bigr)^{-1}
\qquad\text{and}\qquad M(z)=\Gamma_1\gamma(z),   \quad z\in\rho(A_0),
\end{equation}
are called the $\gamma$-field and  Weyl function corresponding to $\Pi$.

\medskip

An important property of the Weyl function  is that $z\in \rho(A_\Theta)\cap \rho(A_0)$ if and only if $0\in \rho(\Theta - M(z)).$
Note that $M\in R^{\rm u}[\K]$ and the following Krein-type formula holds:
\begin{equation}
\label{eq:A.17}
(A_\Theta-z I)^{-1} - (A_0 -z I)^{-1} = \gamma(z) (\Theta - M(z))^{-1} \gamma({\overline
z})^*, \quad z\in \rho(A_0)\cap \rho(A_\Theta).
\end{equation}

\medskip

{\bf A3.~Perturbation determinants.}
Next,  we recall   some basic facts on perturbation determinants treated  in the
framework of boundary triplets (see \cite{MN1} for details).
If  $\wt A',\wt A\in \Ext_A$  and meet the assumptions of Theorem \ref{prop3.3},
then  there is  a  regular boundary triplet for $\{\wt A',\wt A\}$, i.e., $\wt A' =
A_{B'}$,  $\wt A = A_B$, with $B', B \in \cB(\K)$ and by \eqref{eq:A.17}, $B' - B\in
\bS_1$.
Therefore the following expression is well defined on $\rho(\wt A)\cap \rho(A_0)$
$$
\gD^\Pi_{\wt A'/\wt A}(z) = \det(I + (B' -B)(B-M(z))^{-1}, \quad z \in \rho(\wt A)\cap
\rho(A_0).
$$
It is called the {\it perturbation determinant of the pair $\{\wt A',\wt A\}$ with respect
to} $\Pi$, see \cite{MN1}. The following  formula plays the key role in what follows
$$
\big(\gD^\Pi_{\wt A'/\wt A}(z)\big)^{-1}\frac{d}{dz}\gD^\Pi_{\wt A'/\wt A}(z) = \trace\big((\wt
A'-z I)^{-1} - (\wt A-z I)^{-1}\big), \quad z \in \rho(\wt A)\cap \rho(\wt A').
$$
Let now $\{\wt A'',\wt A',\wt A\}$ be a triple of proper  extensions of $A$ such
that each pair is resolvent comparable. Then the following three perturbation determinants are well
defined and

\begin{equation}
\label{eq:A.20}
\gD^\Pi_{\wt A''/\wt A'}(z)  \gD^\Pi_{\wt A'/\wt A}(z) = \gD^\Pi_{\wt A''/\wt A}(z),
\qquad z \in \rho(\wt A)\cap \rho(\wt A')\cap \rho(\wt A'').  
\end{equation}

\medskip

{\bf A.4~The existence of a spectral shift function.}
Now we are ready to improve  Theorem 3.13  from \cite{MalNei2015}.

\begin{thm}
\label{IV.7}
Let $\{L_1, L_0\}$ be a pair of m-accumulative resolvent comparable operators. Suppose also that
$(\rho(L_0)\cup\s_{\rm c}(L_0)) \cap \C_- \not= \varnothing$. Then 

\item[\;\;\rm (i)]
There is a symmetric operator $A,$ $n_+(A) = n_-(A)$, such that $L_1, L_0\in \Ext_A$.
Moreover, there is  a boundary triplet $\Pi\,=\,\{\K, \gG_0, \gG_1\}$ for $A^*$,  regular
for the pair $\{L_1, L_0\}$.

\item[\;\;\rm (ii)]
For each (not necessarily regular) boundary triplet $\Pi$, 
the perturbation determinant $\gD^\Pi_{L_1/L_0}$ admits a (non-unique)
representation
  \begin{equation}
  \label{4.42}
\gD^\Pi_{L_1/L_0}(z) = c\; \exp\left(
\int_\R\left((t-z)^{-1} -t(1+t^2)^{-1}\right)\bs\o_{\rm a}(t)dt
\right), \qquad z \in \C_+.
   \end{equation}
with a \emph{complex function} $\bs\o_{\rm a}\in L^1(\R,(1+t^2)^{-1})$ and  a
constant $c \in \C$.

\item[\;\;\rm (iii)]
If  $\bs\o_a$ satisfies   \eqref{4.42}, then $\bs\o_{\rm a}\in \mathscr S\{L_1,L_0\}$, i.e.,  the following 
formula holds:
\begin{equation}
\label{4.43}
\trace\left((L_1 - z I)^{-1} - (L_0 - z I)^{-1}\right) = - 
\int_\R (t-z)^{-2}\bs\o_{\rm a}(t)\,dt, \quad z \in \C_+.
\end{equation}

\item[\;\;\rm (iv)]
If   $L_0$  (respectively, $L_1$) is self-adjoint, then there is $\bs\o_{\rm a}\in
\mathscr S\{L_1,L_0\}$ satisfying $\IM(\bs\o_{\rm a}(t)) \le 0$ $($respectively, $\IM(\bs\o_{\rm a}(t))
\ge 0)$ a.e. on $\R$.

\item[\;\;\rm (v)]
If $\bs\xi_{\rm a} \in  \mathscr S\{L_1,L_0\}$ is \emph{real-valued}, then for almost all $x\in\R$,
\begin{equation}
\label{5.44}
\bs\xi_{\rm a}(x)  = \frac1\pi \lim_{y\downarrow
0}\IM(\log(\gD^\Pi_{L_1/L_0}(x + iy))).
\end{equation}
\end{thm}

\Pf (i). We put
\begin{equation}
\label{1.14Intro}
Af= Lf = L_0f, \quad f \in \dom(A)
 = \dom(L_1) \cap \dom(L_1^*) \cap \dom(L_0) \cap \dom(L_0^*).
\end{equation}
Clearly, $A$ is a closed symmetric 
operator, $L_1$ and $L_0$ are its proper extensions. Moreover, the assumption
$(\rho(L_0)\cup \sigma_{\rm c}(L_0)) \cap \C_- \not= \varnothing$
implies that $n_+(A) = n_-(A)$. 
Assume for simplicity that $A$ is densely defined.
Since  $(\rho(L_0)\cup \sigma_c(L_0)) \cap \C_- \not= \varnothing$ and $\rho(L_0)\cap
\rho(L_1)\supset \C_+$, the extensions   $L_0$ and $L_1$ meet the hypotheses of Theorem
\ref{prop3.3}. Thus, there exists  a  boundary triplet $\Pi = \{\cH, \gG_0, \gG_1\}$
for $A^*$ regular for the pair  $\{L_1,L_0\}$, i.e. such that
$$
L_1  = A_{B'}= A^*\big|\Ker(\gG_1 - B'\gG_0) \quad\text{and}\quad  L_0  =
A_B = A^*\big|\Ker(\gG_1 - B\gG_0)
$$
where the boundary operators $B',B$ are in $\cB(\K)$ (see \eqref{2.2}).

(ii) By Lemma \ref{prop2.1} (ii), $B'$ and $B$ are accumulative because so are
$A_{B'}$ and $A_{B}$. Besides,   $B'-B\in \bS_1$ because 
$(L_1 -\ri I)^{-1} - (L_0 -\ri I)^{-1}\in \bS_1$.

We set $B''\df\re B'+ \ri\im B$ and $L_2 \df  A^*\big|\Ker(\gG_1 - B''\gG_0)=
A_{B''}$. Clearly, $L_2 \in \Ext_A$ and
by Lemma \ref{prop2.1} (ii),  it is $m$-accumulative because so is $B''$. Since $B'
- B'' =\lb\ri(\im B' -\im B) \in \bS_1$,  
 the perturbation determinant  
 \begin{eqnarray*}
\gD^\Pi_{L_1/L_2}(z) 
= \det(I + \ri(\im B' -\im B)(B'' - M(z))^{-1}), \qquad z \in \C_+,
  \end{eqnarray*}
is well defined. Here $M$ is the  Weyl function corresponding to $\Pi$.  Since
$\im B' -\im B \le -\im B = -\im B''$, Lemma \ref{IV.3a} (ii) guarantees a representation
\begin{equation}
\label{4.44}
\gD^\Pi_{L_1/L_2}(z) = \varkappa \exp\Big(
\int_{\Bbb R}\left((t-z)^{-1} -
t(1+t^2)^{-1}\right)\eta(t)dt\Big), \qquad z \in \C_+,
\end{equation}
for  {\emph {a real-valued}} function $\eta$ in $L^2(\R,(1+t^2)^{-1})$ and
$\varkappa \in \T$.

Next, since  $B'' - B =\re B'-\re B \in \bS_1$, the determinant
$\gD^\Pi_{L_2/L_0}$ is well defined. Since $\re B' -\re B = (\re B' -\re B)^*$,
by  Lemma \ref{IV.2} (ii), we obtain a representation
   \begin{align}
   \label{4.45}
\lefteqn{ \gD^\Pi_{L_2/L_0}(z)
= \det(I + (\re B'-\re B)(B - M(z))^{-1}) } \nonumber \\
& &  =  k\exp\Big(
\int_{\Bbb R}\Big((t-z)^{-1} -
t(1+t^2)^{-1}\Big)\xi(t)dt
\Big), \qquad z \in \C_+,\
  \end{align}
with  a {\emph {real function}} $\xi$ in $L^2(\R,(1+t^2)^{-1}dt)$ and  
$k> 0$. Combining \eqref{4.44} with \eqref{4.45}, using 
the chain rule, and setting  $\bs\o_{\rm a} \df \xi + i\eta$,  we  arrive at 
\eqref{4.42} with $c = k\varkappa$.

(iii)   Formula \eqref{4.43} is an immediate  consequence of \eqref{4.42} and
\eqref{eq:A.20}.

(iv)  Let   $L_0 = L_0^*$. Then $B =B^*$ and $B'' =\re B' = (B'')^*$. Hence $L_2 =
A_{B''}= L_2^*$ and $\im B' - \im B = \im B'\le\0$.  By Lemma \ref{IV.3a}\,(i) with $V_+
= |\im B'|$ and $B = B'$, we have
$$
\gD^\Pi_{L_1/L_2}(z) = \Big(\gD^\Pi_{L_2/L_1}(z)\Big)^{-1} =
\overline{\vk}\exp
\Big(-\int_{\Bbb R}\left((t-z)^{-1}
-t(1+t^2)^{-1}\right)\eta_+(t)dt\Big),
\quad z \in \C_+,
$$
where $\eta_+(t)\ge 0$ a.e. on $\R$. Combining this identity with \eqref{4.45},
applying the chain rule  for determinants,  and setting $\bs\go_{\rm a}(t)\df \xi(t) -\ri\eta_+(t)$,
we arrive at \eqref{4.42} with $\IM(\bs\go_{\rm a}(t))\le 0$, $t \in \R$, and $\varkappa \in \T$.
The case $L_1 = L_1^*$ is similar.

(v)  Substituting $\bs\xi_{\rm a}$ in  \eqref{4.42} in place of $\bs\o_{\rm a}$ and taking
the imaginary part of the logarithms of both sides, we obtain
$$
\IM(\log(\gD^\Pi_{L_1/L_0}(x + \ri y))) = 
\int_{\R}y((t-x)^2 + y^2)^{-1}\bs\xi_{\rm a}(t)dt  + \im c.
$$
Applying the Fatou theorem  as $y\downarrow 0$   we arrive at
\eqref{5.44}. $\bl$

\medskip

{\bf Remark.}
{\bf 1.}  Note that definition \eqref{eq:A.20}  is a special case of a more general
definition  introduced in \cite{MN1} which allows us to treat  the case of pairs $\{A_B,
A_{B'}\} (\subset \Ext_A)$ with unbounded closed boundary operators $B$, $B'$ satisfying the following conditions:

(i) $\dom B=\dom B'$;

(ii) $(B'-B)\bigl(B-M(z)\bigr)^{-1}\in \bs S_1 \quad \text{for}\quad z\in\rho(A_B)\cap
\rho(A_0)$.

In this case the operators $A_B$ and $A_{B'}$  are resolvent comparable.

This definition allows us to extend Lemmata  \ref{IV.3a} and \ref{IV.2}  to the case of
unbounded m-accumulative operators $B$ (in Lemma \ref{IV.3a} it is required  in addition to this
that $\dom B = \dom B^*$).

Note also that the choice of $A=\0$ in \eqref{1.14Intro}  allows us to treat the case of trace class
additive perturbations in the framework of the above definition. Indeed, in this case  a
boundary triplet $\Pi=\{\K,\gG_0,\gG_1\}$ can be chosen so that $M(z)=z$,
$A_B =B,$ and $A_{B'}=B'$,   
and condition (i) implies that $A_{B'} = A_B + V$, where   $V = \clos(B'-B)$. Now
definition \eqref{eq:A.20} coincides with the classical definition
$$
\Delta_{B'/B}(z) = \det\bigl(I + (B'-B)(B-z)^{-1}\bigr), \quad z\in \rho(B),
$$
of the perturbation determinant for the pair $\{B',B\}$, $B'=B+V,$  (see
\cite{Kr}-\cite{Kr87}, \cite{Ya}).

{\bf 2.} In this case ($A=\0$, $M(z)=z$, $A_B =B,$ $A_{B'}=B'$) formula \eqref{4.4a}
takes the form
  \begin{equation}
   \label{4.4a_new}
\!\det\left(I\!+\!(B'-B)(B - z)^{-1}\right)\!=\!c\, \exp\left(\frac 1\pi \int_\R
\left(\frac{1}{t-z} -\frac{t}{1+t^2}\right)\bs\xi(t)dt\right), \!\quad z \in \C_+.  \end{equation}
If in addition to this $V=V^*\in \bs S_1$,  it can be shown (see the proof of Lemma 5.3 of
\cite{MalNei2015}) that $\bs\xi\in L^1(\R)$
 and formula \eqref{4.4a_new} is reduced to
$\Delta_{B'/B}(z) = \exp\left(\frac 1\pi \int_\R (t-z)^{-1}\bs\xi(t)dt\right)$.
 So, if $B = B^*$,  Lemma \ref{IV.2} coincides with the Krein classical result proved originally in
\cite{Kr} by using the step method (see also \cite{Ya}).

{\bf 3.} Note that in the case of additive trace class perturbations the proof of
Theorem \ref{IV.7} is simplified and the technique of boundary triplets is unnecessary.

\medskip

Note in conclusion that several papers are devoted to trace formulae for pairs
$\{-\Delta, \Delta+ q\}$, where $\Delta$ is the Laplacian  operator either in $L^2(\Bbb
R^n)$, or on a discrete lattice and $q$ is a complex-valued potential (see e.g.
\cite{KoLa} and references therein).


%
%
%
%
%
\end{appendix}

\medskip

\footnotesize
\noindent
\begin{tabular}{p{4.6cm}p{4.5cm}p{4.6cm}}
M.M. Malamud  &H. Neidhardt&V.V. Peller \\
Institute of Applied & Institut f\"ur Angewandte & Department of Mathematics\\
Mathematics and Mechanics& Analysis und Stochastik&Michigan State University  \\
NAS of Ukraine &Mohrenstr. 39& East Lansing, Michigan 48824\\
Slavyansk&D-10117 Berlin&USA\\
Ukraine&Germany
\end{tabular}

\end{document}